\documentclass[11pt]{amsart}
\usepackage{amssymb,amsmath,amsfonts,amscd,euscript}

\newcommand{\nc}{\newcommand}

\nc{\slt}{\mathfrak{sl}_2}
\nc{\suth}{\widehat{\mathfrak{su}}(2)}
\nc{\gl}{\mathfrak{gl}}
\nc{\GL}{\mathfrak{GL}}
\nc{\g}{\mathfrak{g}}
\nc{\h}{\mathfrak{h}}
\nc{\la}{\lambda}
\nc{\slth}{\widehat{\slt}}
\nc{\C}{\mathbb C }
\nc{\Z}{\mathbb Z }
\nc{\N}{\mathbb N }
\nc{\R}{\mathbb R }
\nc{\Q}{\mathbb Q }
\nc{\al}{\alpha }
\nc{\be}{\beta}
\nc{\ve}{\varepsilon}
\nc{\ch}{{\mathop {\rm ch}}}
\nc{\Tr}{{\mathop {\rm Tr}\,}}
\nc{\Id}{{\mathop {\rm Id}}}
\nc{\U}{{\mathop {\rm U}}}
\nc{\bra}{\langle}
\nc{\ket}{\rangle}
\nc{\x}{{\bf x}}
\nc{\pa}{\partial}
\nc{\ld}{\ldots}
\nc{\cd}{\cdots}
\nc{\hk}{\hookrightarrow}
\nc{\A}{\mathfrak A}
\nc{\qb}[2]{\genfrac{(}{)}{0pt}{}{#1}{#2}_q}
\nc{\n}{\mathfrak{n}}
\nc{\un}{\mathfrak{u}}
\nc{\T}{\otimes}
\nc{\bv}{{\bf v}}
\nc{\bu}{{\bf u}}
\nc{\bs}{{\bf s}}
\nc{\bw}{{\bf w}}
\nc{\bin}[2]{{\genfrac{(}{)}{0pt}{0}{#1}{#2}}}
\nc{\fac}[1]{(#1)_q!}
\nc{\wt}{\widetilde}

\nc{\bn}{{\bf n}}
\nc{\bc}{{\bf c}}
\nc{\bd}{{\bf d}}
\nc{\bm}{{\bf m}}
\nc{\bl}{{\bf \lambda}}

\newtheorem{theorem}{Theorem}[section]
\newtheorem{lem}{Lemma}[section]
\newtheorem{prop}{Proposition}[section]
\newtheorem{cor}{Corollary}[section]

\newtheorem{conj}{Conjecture}[section]

\title[$2D$ current algebras]
{Two dimensional current algebras and affine fusion product}
\author{B. Feigin, E. Feigin}
\address{BF: Landau institute for Theoretical Physics,\newline
Russia, Moscow region, Chernogolovka, 142432, prospect Ak. Semenova, 1a
 \newline and \newline
Independent University of Moscow,\newline Russia, Moscow, 119002,
Bol'shoi Vlas'evski per., 11}
\email{feigin@mccme.ru}
\address{EF:
Tamm Theory Division, Lebedev Physics Institute,\newline 
Russia, Moscow, 119991,
Leninski pr., 53\newline and\newline
Independent University of Moscow,\newline Russia, Moscow, 119002,
Bol'shoi Vlas'evski per., 11}
\email{evgfeig@mccme.ru}

\begin{document}
\begin{abstract}
In this paper we study a family of commutative algebras generated by two
infinite sets of generators. These algebras are pa\-ra\-me\-trized by Young diagrams. 
We explain a connection of these algebras with the fusion product of 
integrable irreducible representations of the affine $\slt$ Lie algebra.
As an application we derive a fermionic formula for the character
of the affine fusion product of two modules. These fusion products can be
considered as a simplest example of the double affine Demazure modules. 
\end{abstract}
\maketitle

\section*{Introduction}
The main goal of this paper is to derive a formula for the character
of the fusion product of two integrable irreducible representations of
the affine Kac-Moody Lie algebra $\slth$. We first briefly recall the 
definition.

The notion of the fusion product of cyclic representations $V_1,\ldots, V_n$ 
of the Lie algebra  $\g$ was introduced in \cite{FL}. This
object is a cyclic graded representation of the current
algebra $\g\T\C[u]$. 
One starts with a tensor product 
of evaluation representations $V_1(z_1)\T\ldots\T V_n(z_n)$, where $z_i$
are pairwise distinct complex numbers. 
Introduce a filtration  
\begin{equation}
\label{deffus}
F_l=\mathrm{span}
\{(g_1\T u^{i_1}\ldots g_k\T u^{i_k})\cdot v_1\T\ldots\T v_n, \ 
\sum i_\al\le l\},
\end{equation} 
where $v_i$ are cyclic vectors of $V_i$. Then the fusion product is an adjoint 
graded space 
$$V_1(z_1)*\ldots * V_n(z_n)=F_0\oplus\bigoplus_{l>0} F_l/F_{l-1}.$$
An important property is that in some cases $V_1(z_1)*\ldots * V_n(z_n)$ is 
independent on $\{z_i\}$ as $\g\T\C[u]$-module. This is always true for
\begin{itemize}
\item $n=2$ and arbitrary cyclic modules (obvious);
\item $\g=\slt$ and $n$ finite-dimensional modules 
(see \cite{FL}, \cite{FF1});
\item $\g={\mathfrak{sl}_n}$ and irreducible representations with
special highest weights (see \cite{CL},\cite{FKL}, \cite{Ked});
\item simple Lie algebra $\g$ and special highest weight irreducible 
representations (see \cite{FoL}). 
\end{itemize}
In these cases we omit numbers
$z_i$ and denote the corresponding fusion product simply by 
$V_1*\ldots * V_n$.  

Now let $\g$ be an affine Kac-Moody Lie algebra and
$V_i$ be integrable irreducible representations. 
This situation for
$\g=\slth$ and $n=2$ was studied in \cite{FFJMT} in order to derive
some results about monomial bases for vertex operators in minimal models. 
In particular a bosonic formula for the character of fusion product of
level $1$ and level $k$ representations was obtained. Here we consider
two arbitrary level representations and derive a fermionic formula for the 
corresponding fusion product. We briefly describe our approach here.

Let $\g=\slth=\slt\T\C[t,t^{-1}]\oplus\C K\oplus \C d$, where $K$ is
a central element and $d$ is a degree operator. For $k\in\N$ and 
$0\le i\le k$ let $L_{i,k}$ be the corresponding irreducible integrable
representation of $\slth$ with highest weight vector $v_{i,k}$ such that
$h_0v_{i,k}=iv_{i,k}$, $Kv_{i,k}=kv_{i,k}$, $dv_{i,k}=0$,
where $h$ is a standard generator of the Cartan subalgebra and for $x\in\slt$ 
we set $x_i=x\T t^i$. In what follows we use the notation $e,h,f$ for standard
basis of $\slt$.
There exist three different gradings on a fusion product 
$L_{i_1,k_1}*L_{i_2,k_2}$: $\deg_z$ by an operator $h_0$, $\deg_q$ by an
operator $d$ and $\deg_u$ coming from the fusion filtration $(\ref{deffus})$.
This defines a character 
$\ch_{z,u,q} L_{i_1,k_1}* L_{i_2,k_2}$. To find this character we consider
a principal subspace $W_{i,k}^N\hk L_{i,k}$ which is generated from 
the extremal vector $v_{i,k}^N$ by the action of operators $e_i$, $i\in\Z$
(see \cite{FS}, \cite{FF2}).
The weight of $v_{i,k}^N$ is a weight of $v_{i,k}$ shifted by the $N$-th
power of the translation element from the Weyl group of $\slth$. An important
thing is that $L_{i,k}$ is a limit of $W_{i,k}^N$ while $N\to\infty$. This
gives
$$L_{i_1,k_1}* L_{i_2,k_2}=\lim_{N\to\infty} W_{i_1,k_1}^N* W_{i_2,k_2}^N.$$
So we only need to find the character of $W_{i_1,k_1}^N* W_{i_2,k_2}^N$.
We note that this space is a cyclic representation of an abelian Lie algebra 
with  a basis $e_i, e_i\T u$, $i\in\Z$ with cyclic vector 
$v_{i_1,k_1}\T v_{i_2,k_2}$. We give a precise description of 
$W_{i_1,k_1}^N* W_{i_2,k_2}^N$ in terms of generators and relations. 

Let  $\la_0\ge\la_1\ge\cd \ge \la_s>0$ be a set of positive integers.
We consider a corresponding partition
$\bl=\{(i,j),\ i,j\in\Z, i,j\ge 0, i\le\la_j\}$ and define an algebra
\begin{equation}
A_\bl=\C[a_0,a_{-1},\ld; b_0,b_{-1},\ld]/\bra a(z)^ib(z)^j, (i,j)\notin\bl\ket,
\end{equation}
where $a(z)=\sum_{k\ge 0} a_{-k} z^k$, $b(z)=\sum_{l\ge 0} b_{-l} z^l$.
This algebra is two-\-dimen\-si\-onal generalization of an algebra 
$\C[e_i]_{i\le 0}/e(z)^{k+1}$ which plays an important role in the 
representation theory of $\slth$ (see for example \cite{FS}). 
Recall that one of the most useful tools of study of $A_k$ and similar
algebras is a vertex operator realizations (see for example \cite{FK},
\cite{FF3}, \cite{FST}). We also use this technique in our situation. Namely
for $\bl$ such that  
$\la_{i-1}-\la_i\le \la_i-\la_{i+1}$  we embed $A_\bl$ into multi-dimensional
lattice vertex operator algebra and compute the character of $A_\bl$. 
In addition we show that there exists $\bl$ (depending on $k_1$ and $k_2$) such that
up to a certain identification of generators $e_i$, $e_i\T u$ and 
$a_j, b_j$  the space $W_{i_1,k_1}^N* W_{i_2,k_2}^N$ is isomorphic to
the quotient of  
$A_\bl$ by some ideal. 
This ideal depends on $i_1$ and $i_2$ and is generated by certain coefficients of the series
$a(z)^i$ and $b(z)^j$.
Using the realization above we derive a fermionic (Gordon type)  formula for 
the the character of the fusion
product of principal subspaces and therefore (as a limit) the character of
the fusion product of two integrable irreducible representations. For
fermionic formulas in the case of finite-dimensional algebras see for
example \cite{FF1, FJKLM1, FJKLM2, AK}.

We also expect the existence of the bosonic (alternating sign) formula for 
the fusion product 
of integrable representations  similar to one given
in \cite{FFJMT} in the $\slth$ case with one of the representations of 
the level $1$. Let us briefly explain the importance of such formula. Let
$\g$ be an affine Kac-Moody Lie algebra, $L_\mu, L_\nu$ be its integrable
highest weight representations. Then one has a decomposition of the 
tensor product
$$L_\mu\T L_\nu=\bigoplus C_{\mu,\nu}^\pi L_\pi$$
into the direct sum of integrable highest weight modules $L_\pi$, where
$C_{\mu,\nu}^\pi$ are the spaces of multiplicities (see \cite{K1}). 
The fusion filtration $(\ref{deffus})$ defines a filtration on 
$C_{\mu,\nu}^\pi$, because all $F_l$ in $(\ref{deffus})$ are 
representations of $\g=\g\T 1\hk\g\T \C[t]$. We note that spaces of
multiplicities appear in the coset conformal field theories (see \cite{DMS}).
In \cite{FFJMT} the character of the filtered space  $C_{\mu,\nu}^\pi$
in some particular case
was given in the bosonic form and the resulting formula was applied for the 
study of vertex operators in Virasoro minimal models. We hope that 
in the general case (at least for $\g=\slth$) it is also possible to write
an analogous bosonic formula and establish connection with the 
corresponding coset conformal field theory.

In the end of the introduction we comment about the connection of the affine 
fusion  product with two-dimensional affine Demazure modules. Recall that 
in some special cases affine Demazure modules 
(see for example \cite{Kum})
can be identified with 
the fusion product of finite-dimensional irreducible representations
of the corresponding simple algebra (see \cite{FF2}, \cite{CL}, \cite{FKL},
\cite{FoL}). In addition some
integrable irreducible representations (in particular the basic one)
can be realized as an inductive limit of these fusion products 
(see \cite{FF2, FoL}).
All these statements allow to consider affine fusion products of two 
representations as  a simplest example of double affine Demazure modules.
To construct another examples one needs to "fuse" $N$ irreducible 
representations
for arbitrary $N$. In particular it seems to be important to prove the 
independence of the corresponding fusion products on the evaluation parameters.
We hope to return to these questions elsewhere.

Our paper is organized in the following way.

In Section $1$ we recall the definition and collect main properties of the 
lattice vertex operator algebras and principal subspaces.

In Section $2$  we study the family of commutative algebras labeled by Young 
diagram and construct their vertex operator realization.

In Section $3$ we apply the results of Section $2$ to the computation of the 
character of the fusion product of two irreducible representations.

{\it Acknowledgements.}\\ 
Research of BF is partially supported by RFBR Grants 04-01-00303 and
05-01-01007, INTAS 03-51-3350 and NSh 2044.2003.2.
Research of EF is partially supported by the 
RFBR Grant 06-01-00037 and LSS 4401.2006.2.

\section{Lattice vertex operator algebras}
In this section we recall main properties of lattice vertex
operator algebras (VOA for short) and derive some statements about
principal subspaces.
The main references are \cite{K2}, \cite{BF},
\cite{dong}.

Let $L$ be a lattice of finite rank equipped with a symmetric bilinear form
$(\cdot,\cdot): L\times L\to \Z$ such that $(\al,\al)>0$ for  all $\al\in 
L\setminus \{0\}$. Let $\h=L\T_{\Z} \C$. The form  $(\cdot,\cdot)$ induces
a bilinear form on $\h$, for which we use the same notation. Let 
$$\widehat\h=\h\T\C[t,t^{-1}]\oplus \C K$$
be the corresponding multi-dimensional Heisenberg algebra with the bracket
$$[\al\T t^i, \be\T t^j]=i\delta_{i,-j}(\al,\be) K,\  
[K,\al\T t^i]=0,\ \al,\be\in\h.$$
For $\al\in\h$ define the 
Fock representation $\pi_\al$ generated by a vector $|\al\ket$ such that
$$(\be\T t^n) |\al\ket=0,\ n>0; \qquad (\be\T 1) |\al\ket=(\be,\al) |\al\ket;
\qquad
K|\al\ket=|\al\ket.$$

We now define a VOA $V_L$ associated with $L$. We deal only with an even
case, i.e.  $(\al,\be)\in 2\Z$ for all $\al,\be\in L$ 
(in the general case the construction leads to the 
so called super VOA). As a
vector space 
$$V_L\simeq \bigoplus_{\al\in L} \pi_\al.$$
The $q$-degree on $V_L$ is defined by 
\begin{equation}
\label{defqdeg}
\deg_q |\al\ket= \frac{(\al,\al)}{2},\quad \deg_q (\be\T t^n)=-n.
\end{equation}
The main ingredient of the VOA structure 
on $V_L$ are bosonic vertex 
operators
$\Gamma_\al(z)$ which correspond to highest weight vectors $|\al\ket$.
One sets
\begin{equation}
\label{defvo}
\Gamma_\al(z)=S_\al z^{\al\T 1} \exp(-\sum_{n<0} \frac{\al\T t^n}{n} z^{-n})
\exp(-\sum_{n>0} \frac{\al\T t^n}{n} z^{-n}),
\end{equation}
where $z^{\al\T 1}$ acts on $\pi_\be$ by $z^{(\al,\be)}$ and 
the operator $S_\al$ is defined by
$$S_\al |\be\ket=c_{\al,\be} |\al+\be\ket;\quad
[S_\al,\be\T t^n]=0, \ \al,\be\in\h,$$
where $c_{\al,\be}$ are some nonvanishing constants.
The Fourier decomposition is given by
$$\Gamma_\al(z)=\sum_{n\in\Z} \Gamma_\al(n) z^{-n-(\al,\al)/2}.$$
In particular,
\begin{equation}
\label{htoh}
\Gamma_\al(-(\al,\al)/2-(\al,\be))|\be\ket=c_{\al,\be}|\al+\be\ket.
\end{equation}
One of the main properties of vertex operators is the following
commutation relation:
\begin{equation}
\label{vertop}
[\al\T t^n,\Gamma_\be(z)]=(\al,\be)z^n \Gamma_\be(z).
\end{equation}
Another important formula describes the product of two vertex operators
\begin{multline}
\Gamma_\al(z)\Gamma_\be(w)=(z-w)^{(\al,\be)} S_\al S_\be 
z^{(\al+\be)\T 1}\times\\ 
\exp(-(\sum_{n<0} \frac{\al\T t^n}{n} z^{-n} + \frac{\be\T t^n}{n} w^{-n}))
\exp(-(\sum_{n>0} \frac{\al\T t^n}{n} z^{-n} +  \frac{\be\T t^n}{n} w^{-n})).
\end{multline}
This leads to the proposition:

\begin{prop}
\begin{equation}
\label{vorel}
(\Gamma_\al(z))^{(k)}(\Gamma_\be(z))^{(l)}=0 \text{ if } k+l<(\al,\be),
\end{equation}
where the superscript $(k)$ denotes the $k$-th derivative of the corresponding
series. In addition if $(\al,\be)=0$ then 
$$\Gamma_\al(z)\Gamma_\be(z) \text{ is proportional to } \Gamma_{\al+\be}(z).$$
\end{prop}

We now recall some basic facts about the representation theory of $V_L$.
This VOA is known to be rational, i.e. every $V_L$-module is completely
reducible. The number of irreducible representations is finite. These 
representations 
are labeled by the elements of $L'/L$, where $L'$ is a dual lattice
\begin{equation}
\label{duallat}
L'=\{\be\in L\T_\Z \Q:\ (\al,\be)\in\Z\ \forall \al\in L\}.
\end{equation}
Namely for $\gamma\in L'/L$ define
$$V^\gamma_L=\bigoplus_{\be\in L+\gamma} \pi_\be.$$
For example vertex operators $\Gamma_\al(z)$ act on each $V^\gamma_L$ via
the formula $(\ref{defvo})$ (because $(\al,\be+\gamma)\in\Z$ for all 
$\be\in L$). The $q$-degree on $V^\gamma_L$ is defined as in ($\ref{defqdeg}$).

In what follows we fix a set $\al_1,\ldots,\al_N$
of linearly independent vectors generating 
the lattice $L$. We denote the nondegenerate matrix of scalar products by
$M$: $m_{ij}=(\al_i,\al_j)$ and assume that $m_{ij}\in \N$ . 

\begin{lem}
For any vector $\bv\in \Z_{\ge 0}^N$ there exist $\gamma_\bv\in L'/L$  and
$\be\bv\in L+ \gamma_\bv$ such that
\begin{equation}
\label{begcond}
\Gamma_{\al_i} (-n) |\be_\bv\ket=0 \text{ for } n<v_i, 1\le i\le N; \quad
\Gamma_{\al_i} (-v_i) |\be_\bv\ket=c_i |\be_\bv+\al_i\ket, 
\end{equation}
with some nontrivial constants $c_i$.
\end{lem}
\begin{proof}
We only need to find $\be_\bv$ such that 
\begin{equation}
\label{degcond}
\deg_q |\be_\bv+\al_i\ket - \deg_q |\be_\bv\ket=v_i,\ 1\le i\le N.
\end{equation}

Note that $\be\in L\T_\Z \Q$. So $\be_\bv$ is a rational linear combination of
$\al_i$ and we can consider $\be_\bv$ as a vector in $\Q^N$. Then
$(\ref{degcond})$ is equivalent to
\begin{equation}
\label{bebv}
\frac{m_{ii}}{2}+ (M\be_\bv)_i=v_i.
\end{equation}
In view of $m_{ii}/2\in\Z$ we obtain that $\be_\bv\in L'$ satisfying 
$(\ref{bebv})$ really exists, 
because $(\ref{duallat})$ can be rewritten as
$$L'=\{\be\in L\T_\Z \Q:\ M\be\in\Z\}.$$
Then $\gamma_\bv$ is defined as a class of $\be_\bv$.
\end{proof}

We now define principal subspaces. For $\bv\in Z_{\ge 0}^N$ consider 
the subspace $W_L(\bv)\hk V^{\gamma_\bv}_L$ generated from the vector
$\be_\bv$ by an action of operators $\Gamma_{\al_i}(-n_i)$ with
$n_i\ge v_i$  ($1\le i\le N$). Our goal is to describe $W_L(\bv)$ (in particular we want to find its character).  
We first realize this subspace as a quotient of a polynomial
algebra. Namely define $W'_L(\bv)$ as a quotient of 
$\C[a_i(-n)]_{\substack{1\le i\le N\\ n\ge v_i}}$ by relations
$$a_i(z)^{(k)}a_j(z)^{(l)},\ k+l<m_{ij},$$
where $a_i(z)=\sum_{n\ge v_i} z^n a_i(-n)$.
$W_L'(\bv)$ is generated by coefficients 
$$a_i(-n),\quad 1\le i\le N, n\ge v_i.$$
We note that 
$W_L'(\bv)=\bigoplus_{\bn\in\Z_{\ge 0}^N} W_{L,\bn}'(\bv)$, 
where $W_{L,\bn}'(\bv)$ is a subspace spanned by monomials in $a_i(k)$ such that
the number of factors of the type $a_{i_0}(k)$ with fixed $i_0$ is exactly
$n_{i_0}$. The character of $W_{L,\bn}'(\bv)$ is naturally defined by
$\deg_q a_i(k)=-k.$

\begin{lem}
\label{charlem}
\begin{equation}
\label{charact}
\ch_q W_{L,\bn}'(\bv)=
\frac{q^{\bn M\bn/2+\sum_{i=1}^N n_i(v_i-m_{ii}/2)}}{(q)_\bn},
\end{equation}
where $(q)_\bn=\prod_{j=1}^N (q)_{n_j}$, $(q)_n=\prod_{j=1}^n (1-q^j)$.
\end{lem}
\begin{proof}
We use the dual space approach (see \cite{FS, FJKLM1}). 
For $\theta\in (W_{L,\bn}'(\bv))^*$ define
a polynomial 
$f_\theta\in\C[z_{i,n}]_{\substack{1\le i\le N\\ 1\le n\le n_i}}$ by
$$f_\theta=\theta\left(\prod_{i=1}^N \prod_{n=1}^{n_i} a_i(z_{i,n})\right).$$
From the exact form of the relations in $W_L'(\bv)$ one can see that the space
$$\{f_\theta:\ \theta\in (W_{L,\bn}'(\bv))^*\}$$ 
coincides with the space of polynomials
which are divisible by
\begin{multline}
\label{pol}
\prod_{i=1}^N\prod_{n=1}^{n_i} z_{i,n}^{v_i}
\prod_{i=1}^N  \prod_{1\le n< m\le n_i} (z_{i,n}-z_{i,m})^{m_{ii}}
\times\\
  \prod_{1\le i<j\le N} \prod_{\substack{1\le n\le n_i\\ 
       1\le m\le n_j}} (z(i,n)-z(j,m))^{m_{ij}}.
\end{multline}
The character of such polynomials  coincides with the right hand side of
$(\ref{charact})$.
\end{proof}

In the next proposition we show that spaces $W_L(\bv)$ and $W_L'(\bv)$ are
isomorphic.

\begin{prop}
\label{charvo}
The map $\be_\bv\mapsto 1$, $\Gamma_{\la_i}(n)\mapsto a_i(n)$ 
induces the isomorphism
$$W_L(\bv)\simeq W_L'(\bv).$$ 
In particular for any $\bn=(n_1,\ldots,n_N)\in\Z_{\ge 0}^N$
\begin{equation}
\label{chvo}
\ch_q (W_L(\bv)\cap \pi_{\be_\bv+n_1\al_1+\ldots +n_N\al_N})=
\frac{q^{\bn M\bn/2+\sum_{i=1}^N n_i(v_i-m_{ii}/2)}}{(q)_\bn} 
q^{\deg_q \be_\bv}.
\end{equation}
\end{prop}
\begin{proof}
Because of $(\ref{vorel})$ it suffices to prove the equality
$(\ref{chvo})$.

We note that  
$$|\be_\bv+n_1\al_1+\ldots +n_N\al_N\ket \in W_L(\bv)$$
(see $(\ref{htoh})$).
In addition 
\begin{multline*}
\deg_q|\be_\bv+\sum_{i=1}^N n_i\al_i\ket- \deg_q|\be_\bv\ket=
\frac{1}{2}((\be_\bv+\bn) M (\be_\bv+\bn)- \be_\bv M \be_\bv=\\
\frac{1}{2}(\bn M \bn +2 \bn M \be_\bv)=
\frac{1}{2}\bn M \bn + \sum_{i=1}^N n_i (v_i-\frac{m_{ii}}{2}),
\end{multline*}
where the last equality is true because of $(\ref{bebv})$.
Therefore the minimal power of $q$ in the left hand side of $(\ref{chvo})$
is equal to 
$$\bn M\bn/2+\sum_{i=1}^N n_i(v_i-m_{ii}/2).$$

For $\theta\in (W_L(\bv)\cap \pi_{\be_\bv +n_1\al_1+\ldots + n_N\al_N})^*$ set
$$f_\theta(z_{i,n})_{\substack{1\le i\le N\\ 1\le n\le n_i}}=
\theta(\prod_{i=1}^N\prod_{n=1}^{n_i} \Gamma_{\al_i}(z_{i,n})).$$
In view of $(\ref{vertop})$ we obtain that 
$W_L(\bv)\cap \pi_{\be_\bv + n_1\al_1+\ldots +n_N\al_N}$ is invariant under 
the action of operators $\al_i\T t^k$, $k\ge 0$, $1\le i\le N$. 
For the action on the dual
space one has
$$(\al_i\T t^k)f_\theta=
(\sum_{j=1}^N m_{ij} \sum_{n=1}^{n_j} z_{j,n}^k) f_\theta.$$
Using the nondegeneracy of $M$ and a fact that polynomials
$\sum_{n=1}^{n_j} z_{j,n}^k$, $k\ge 0$ generates the ring of symmetric 
polynomials  in variables $z_{j,n}$ with fixed $j$ we obtain that the 
character of $W_L(\bv)\cap \pi_{\be_\bv+n_1\al_1+\ldots +n_N\al_N}$ is greater 
than or 
equal to the right hand side of $(\ref{chvo})$. Now using
$(\ref{vorel})$ and Lemma $\ref{charlem}$ we obtain our proposition.
\end{proof}

\section{Commutative algebras}
Let  $\la_0\ge\la_1\ge\cd \ge \la_s>0$ be a set of positive integers.
We consider a corresponding partition
$\bl=\{(i,j),\ i,j\in\Z_{\ge 0}, i\le\la_j\}$ and define an algebra
\begin{equation}
A_\bl=\C[a_0,a_{-1},\ld; b_0,b_{-1},\ld]/\bra a(z)^ib(z)^j, (i,j)\notin\bl\ket,
\end{equation}
where $a(z)=\sum_{k\ge 0} a_{-k} z^k$, $b(z)=\sum_{l\ge 0} b_{-l} z^l$.
This means that our algebra is generated by two sets of variables $a_k$ and
$b_l$ and the ideal of relations is generated by coefficients of series
$a(z)^ib(z)^j$, $(i,j)\notin\bl$.
For example, $a(z)^{\la_0+1}=b(z)^{s+1}=0$.
$A_\bl$ is graded by
$$\deg_z a_n= \deg_z b_n=1,\ \ \deg_u a_n=0, \deg_u b_n=1,\ \  
\deg_q a_n=\deg_q b_n=-n.$$ 
We want to find
the character of $A_\bl$, which is given  by
\begin{multline}
\ch_{z,u,q} A_\bl=\\ \sum_{i,j,k \ge 0} z^i u^j q^k 
\dim \{x\in A_\bl:\ \deg_z x=i, \deg_u x=j, \deg_q x=k\}.
\end{multline}

We first recall the corresponding result in one-dimensional case
(see \cite{FS, FJKLM1}).

\begin{prop}
\label{filtr}
Let $A_k=\C[a_0,a_{-1},\ld]/a(z)^{k+1}$. Then there exists a Gordon filtration
$F_\mu$ of $A_k$ (labeled by Young diagrams $\mu$)
such that the adjoint graded algebra is
generated by coefficients of series $a^{[i]}(z)$, which are images of
powers $a(z)^i$, $1\le i\le k$. 
In addition defining relations in the adjoint graded algebra
are of the form
\begin{equation}
\label{rel}
a^{[i]}(z)^{(l)} a^{[j]}(z)^{(r)}=0 \text{ if } s+r< 2\min(i,j).
\end{equation}
Here the superscript $(l)$ is used for the $l$-th derivative of the
corresponding series. Relations $(\ref{rel})$ gives the following Gordon type
formula:
\begin{equation}
\label{G}
\ch_{z,q} A_k =\sum_{\bn\in \Z_{\ge 0}^k} (zq^{-1})^{|\bn|}
\frac{q^{\bn A\bn/2}}{(q)_\bn},
\end{equation}
where $A_{i,j}=2\min(i,j)$, $|\bn|=\sum_{i=1}^k in_i$ and
\begin{equation}
(q)_\bn=\prod_{i=1}^k (q)_{n_i}, \ (q)_j=\prod_{s=1}^j (1-q^s).
\end{equation}
\end{prop}

The following theorem gives two-dimensional generalization of $(\ref{G})$
for special $\bl$.

\begin{theorem}
\label{mth}
We consider $\la_0\ge \ld\ge \la_s > 0$ with a condition
\begin{equation}
\la_{i-1}-\la_i\le \la_i-\la_{i+1},\ i=1,\ld, s-1.
\end{equation}
Then
\begin{equation}
\label{mf}
\ch_{z,u,q} A_\bl=\sum_{\bn\in \Z^{\la_0}_{\ge 0}, \bm\in\Z_{\ge 0}^{s}}
u^{|\bm|} (zq^{-1})^{|\bn|+|\bm|} 
\frac{q^{\bn A\bn/2 + \bn B\bm + \bm A\bm/2}}
{(q)_\bn (q)_\bm},
\end{equation}
where $A_{i,j}=2\min(i,j)$ and $B_{i,j}=\max(0,i-\la_j)$.
\end{theorem}

We first construct an algebra with the character given exactly by the right
hand side of $(\ref{mf})$. Let $\A_\bl$ be an algebra generated by
coefficients of abelian currents
$$a^{[i]}(z)=\sum_{n\ge 0} a^{[i]}_{-n} z^n,\ 1\le i\le \la_0,\quad
  b^{[j]}(z)=\sum_{n\ge 0} b^{[j]}_{-n} z^n,\ 1\le j\le s.$$
These currents are subject to the only relations
\begin{gather}
\label{f}
a^{[i]}(z)^{(l)} a^{[j]}(z)^{(r)}=0 \text{ if } l+r< 2\min(i,j),\\
\label{s}
b^{[i]}(z)^{(l)} b^{[j]}(z)^{(r)}=0 \text{ if } l+r< 2\min(i,j),\\
\label{t}
a^{[i]}(z)^{(l)} b^{[j]}(z)^{(r)}=0 \text{ if } l+r< i-\la_j.
\end{gather}
We define $3$ gradings on on $\A_\bl$:
\begin{gather*}
\deg_z a^{[i]}_n=\deg_z b^{[i]}_n=i,\\ 
  \deg_u a^{[i]}_n=0, \deg_u b^{[i]}_n=i,\\ 
  \deg_q a^{[i]}_n=\deg_q b^{[i]}_n=-n.
\end{gather*}
This defines the character $\ch_{z,u,q} \A_\bl$ of $\A_\bl$.

\begin{lem}
\label{dual}
$\ch_{z,u,q} \A_\bl$ coincides with the right hand side of $(\ref{mf})$.
\end{lem}
\begin{proof}
Follows from Lemma $\ref{charlem}$.
\end{proof}

We now prove our theorem in two steps comparing left and right hand sides
of $(\ref{mf})$. 

\begin{lem}
For any $\bl$ the left hand side of $(\ref{mf})$ is less than or equal 
to the right hand side.
\end{lem}
\begin{proof}
Consider a filtration $F_\mu(A_{\la_0})\cdot \C[b_i]_{i=0}^\infty$ on $A_\bl$ 
(this filtration comes from the
filtration  of from Proposition $\ref{filtr}$ on the subalgebra of $A_\bl$ generated by
coefficients of $a(z)$).
The adjoint graded algebra $A'_\bl$ is generated by coefficients of series
$a^{[i]}(z)$ (images of $a(z)^i$) and $b(z)$. We now consider a filtration
on $A'_\bl$ coming from the filtration from Proposition $\ref{filtr}$ on the subalgebra of
$A'_\bl$ generated by $b_j$. We denote the adjoint graded algebra by
$A^{gr}_\bl$. This algebra is generated by coefficients of currents
\begin{equation}
\label{gen}
a^{[i]}(z), b^{[j]}(z), \quad i\le \la_0,\ j\le s.
\end{equation}

First note that in view of  $(\ref{rel})$ currents $(\ref{gen})$ are subject
to the relations
\begin{gather}
\label{1}
a^{[i]}(z)^{(l)} a^{[j]}(z)^{(r)}=0 \text{ if } l+r< 2\min(i,j).\\
\label{2}
b^{[i]}(z)^{(l)} b^{[j]}(z)^{(r)}=0 \text{ if } l+r< 2\min(i,j).
\end{gather}
We show that in addition
\begin{equation}
\label{3}
a^{[i]}(z)^{(l)} b^{[j]}(z)^{(r)}=0 \text{ if } l+r< \max(0,i-\la_j).
\end{equation}
Note that it is enough to show that relations
\begin{equation}
\label{simp}
(a(z)^i)^{(l)} (b(z)^j)^{(r)}=0 \text{ if } l+r< \max(0,i-\la_j)
\end{equation}
hold in $A_\la$. We use the induction on $i$. For $i=\la_j+1$ we 
need to show that $a(z)^{\la_j+1} b(z)^j=0$. 
But this is a relation in $A_\bl$. Now suppose that ($\ref{simp}$) is proved
for all $i\le i_0$. We want to show that
\begin{equation}
\label{tp}
(a(z)^{i_0+1})^{(l)} (b(z)^j)^{(r)}=0 \text{ if } l+r< i_0-\la_j+1.
\end{equation}
If $l+r <i_0-\la_j$ then $(\ref{tp})$ holds by
induction assumption, because
\begin{equation}
(a(z)^{i_0}a(z))^{(l)} (b(z)^j)^{(r)}=
\sum_{\gamma=0}^l (a(z)^{i_0})^{(\gamma)} (b(z)^j)^{(r)} x_\gamma(z)
\end{equation}
for some $x_\gamma$. Suppose $l+r=i_0-\la_j$. \\
{\bf Case $1$.} $l\ne 0$. Then for some $x(z)$
$$(a(z)^{i_0+1})^{(l)} (b(z)^j)^{(r)}=
  a(z)^{i_0+1-l} (b(z)^j)^{(r)}x(z).$$
In view  of $i_0+1-l\le i_0$ we can use the induction assumption,
which gives
$$a(z)^{i_0+1-l} (b(z)^j)^{(r)}=0,$$
because $l+r=i_0-\la_j$ and so $r<(i_0+1-l)-\la_j$.\\
{\bf Case $2$.} $l=0$. We need to show that
$a(z)^{i_0+1} (b(z)^j)^{(i_0-\la_j)}=0$.
Note that $a(z)^{i_0+1} b(z)^j=0$ for $i_0\ge \la_j$. Therefore
$$(a(z)^{i_0+1} b(z)^j)^{(i_0-\la_j)}=0.$$
But the following equality holds in $A_\bl$:
$$(a(z)^{i_0+1} b(z)^j)^{(i_0-\la_j)}=a(z)^{i_0+1} (b(z)^j)^{(i_0-\la_j)}.$$
In fact 
$$(a(z)^{i_0+1} b(z)^j)^{(i_0-\la_j)}=\sum_{l=0}^{i_0-\la_j}
\bin{i_0-\la_j}{l} (a(z)^{i_0+1})^{(l)} (b(z)^j)^{(i_0-\la_j-l)}.$$
But if $l\ne 0$ then 
$$(a(z)^{i_0+1})^{(l)} (b(z)^j)^{(i_0-\la_j-l)}=0$$
because of the Case $1$. We thus obtain 
$$a(z)^{i_0+1} (b(z)^j)^{i_0-\la_j}=0.$$
Equalities $(\ref{tp})$ and $(\ref{simp})$ are proved.

We now consider an algebra generated by currents
$(\ref{gen})$ with only relations $(\ref{1})$, $(\ref{2})$, $(\ref{3})$.
Because of Lemma $\ref{dual}$ the character of this algebra is given by 
the right hand side of $(\ref{mf})$. This finishes the proof of our lemma.
\end{proof}

To prove an equality in $(\ref{mf})$ we construct a vertex operator realization of 
an algebra $A_\la$.  
Consider a vector space $\h\simeq\R^N$ 
equipped with a standard scalar product $(\cdot,\cdot)$. We fix 
$N$ such that there exists a set of
linearly independent vectors $p_1,\ld, p_{\la_0}, q_1,\ld,q_s\in\R^N$ with
the scalar products $(p_i,p_j)=2\delta_{i,j}$,
$(q_i,q_j)=2\delta_{i,j}$ and
\begin{equation}
\label{scal}
(p_i,q_j)=1 \text{ if } \la_{j-1}\ge i >\la_j;\quad
(p_i,q_j)=0 \text{ otherwise }.
\end{equation}
(For example, take $N=\la_0+s$ and let $e_1,\ld,e_N$ be some orthonormal basis.
Put
$$q_j=\sqrt{2} e_j, 1\le j\le s;\quad p_i=\frac{1}{\sqrt{2}} e_j+
\sqrt{\frac{3}{2}} e_{s+i},\ j \text{ such that } \la_{j-1}\ge i >\la_j.$$
Then these vectors obviously satisfy (\ref{scal}).)
In what follows we fix a lattice $L$ generated by vectors $p_i,q_j$. 
Let $\Gamma_{p_i}(z)$, $\Gamma_{q_j}(z)$ be corresponding bosonic vertex 
operators. Set
\begin{equation}
\label{tatb}
\tilde a(z)=\sum_{i=1}^{\la_0} \Gamma_{p_i}(z),\quad 
\tilde b(z)=\sum_{j=1}^s \Gamma_{q_j}(z).
\end{equation}

\begin{lem}
\label{grel}
Suppose that 
\begin{equation}
\label{cond}
\la_{i-1}-\la_i\le \la_i-\la_{i+1},\ i=1,\ld, s-1.
\end{equation}
Then $\tilde a(z)^i \tilde b(z)^j=0$ for $(i,j)\notin\bl$.
\end{lem}
\begin{proof}
It suffices to check that
\begin{equation}
\label{van}
\Gamma_{p_{l_1}}(z) \ldots \Gamma_{p_{l_{\la_i+1}}}(z) 
\Gamma_{q_{r_1}}(z) \ldots \Gamma_{q_{r_i}}(z)=0
\end{equation}
for any $l_1,\ldots,l_{\la_i+1}$, $r_1,\ldots,r_i$. 
Note that we can assume that 
$$l_1<\ldots <l_{\la_i+1},\quad r_1 <\ldots <r_i,$$ 
because $(p_l,p_l)=2=(q_r,q_r)$  
and therefore  $\Gamma_{p_l}(z)^2=\Gamma_{q_r}(z)^2=0$. 
In view of  $(\ref{scal})$ we have
$$\Gamma_{p_l}(z)\Gamma_{q_r}(z)=0 \text{ for } \la_{r-1}\ge l >\la_r.$$
So $(\ref{van})$ holds if there exists $k\le \la_r+1$ such that
$\la_{\be_t-1}\ge l_k >\la_{\be_t}$ for some $t\le r$. The number of
such $k$ equals to
$$\sum_{t=1}^r (\la_{\be_t}-\la_{\be_t-1}).$$
Because of the condition $(\ref{cond})$ this sum is greater than or equal to
$\la_0-\la_r$. So we obtain $(\ref{van})$, because the number of factors of the 
type $\Gamma_{p_l}(z)$ is equal to $\la_r+1$.  Lemma is proved.
\end{proof}

This Lemma can be used for the proof of $(\ref{mf})$. But 
in the last section we will need a modification of Theorem $\ref{mth}$.
So we formulate and prove its slight generalization.

Let 
$$
\bc=(c_1,\ldots, c_{\la_0})\in \Z_{\ge 0}^{\la_0},\quad
\bd=(d_1,\ldots, d_s)\in\Z_{\ge 0}^s.$$
We consider an ideal $I_{\la;\bc,\bd}\hk A_\la$ generated by conditions
\begin{gather}
\label{*}
a(z)^i\div z^{c_i+2c_{i-1}+\ldots + ic_1}, \qquad 1\le i\le \la_0,\\
\label{**}
b(z)^j\div z^{d_j+2d_{j-1}+\ldots + jd_1}, \qquad 1\le j\le s.
\end{gather}
This means $I_{\la;\bc,\bd}$ is generated by 
coefficients of $a(z)^i$ 
in front of the powers $z^r$, $0\le r< \sum_{l=1}^i lc_{i+1-l}$  and
by coefficients of $b(z)^j$ 
in front of the powers $z^r$, $0\le r< \sum_{l=1}^j ld_{j+1-l}$.
We define
$$A_{\la;\bc,\bd}=A_\la/I_{\la;\bc,\bd}.$$

\begin{theorem}
Let
\begin{equation}
\la_{i-1}-\la_i\le \la_i-\la_{i+1},\ i=1,\ld, s-1.
\end{equation}
Then
\begin{multline}
\label{gmf}
\ch_{z,u,q} A_{\la;\bc,\bd}=\\
\sum_{\bn\in \Z^{\la_0}_{\ge 0}, \bm\in\Z_{\ge 0}^{s}}
u^{|\bm|} (zq^{-1})^{|\bn|+|\bm|}\times\\ 
\frac{q^{\bn A\bn/2 + \bn B\bm + \bm A\bm/2}}
{(q)_\bn (q)_\bm}
q^{\sum\limits_{1\le i\le j\le\la_0} (j-i+1)n_j c_i +
\sum\limits_{1\le i\le j\le s} (j-i+1)m_jd_i},
\end{multline}
where $A_{i,j}=2\min(i,j)$ and $B_{i,j}=\max(0,i-\la_j)$.
\end{theorem}

\begin{proof}
Recall that $L$ is a lattice generated by vectors $p_i$, $q_j$. 
Consider a vector $\bv\in L$
$$\bv=c_1p_1+(c_1+c_2)p_2 + \ldots + (c_1+\ldots + c_{\la_0})p_{\la_0} +
d_1 q_1 + \ldots + (d_1+\ldots +d_s)q_s$$
and the corresponding principal subspace $W_L(\bv)$.
Denote 
$$\tilde\Gamma_\al(z)=\sum_{n\ge 0} \Gamma_{p_i}(-n) z^n,$$
which differs from usual vertex operators by a certain power of $z$.
Then for 
\begin{equation}
\label{point}
\tilde a(z)=\sum_{i=1}^{\la_0} \tilde \Gamma_{p_i}(z),\quad
  \tilde b(z)=\sum_{j=1}^s \tilde \Gamma_{q_j}(z)
\end{equation}
one gets 
$$\tilde a(z)^i\cdot |\be_\bv\ket\div z^{ic_1+\ldots +c_i},\quad
  \tilde b(z)^j\cdot |\be_\bv\ket\div z^{id_1+\ldots +d_i}.$$
Combining these relations with Lemma $\ref{grel}$ we obtain a
surjection
$A_{\bl;\bc,\bd}\to W_L(\bv)$. We now construct the degeneration of 
$W_L(\bv)$ in order to show that this surjection is an isomorphism.

Consider a family of spaces $W_L(\bv,\ve)$ ($\ve\in\R$, $\ve>0$) 
generated from the vector $|\be_\bv\ket$ by coefficients
of series
\begin{equation}
\label{degen}
\tilde a_\ve(z)=\sum_{i=1}^{\la_0} \ve^i \tilde\Gamma_{p_i}(z),\quad
 \tilde b_\ve(z)=\sum_{j=1}^s \ve ^j \tilde\Gamma_{q_j}(z).
\end{equation}
For any positive $\ve$  a space $W_L(\bv,ve)$ is isomorphic to  $W_L(\bv)$.
Denote the limit $\lim_{\ve\to 0} W_L(\bv,\ve)$ by $W_L(\bv,0)$.
Then 
the $q$-characters of $W_L(\bv,0)$ and $W_L(\bv)$  coincide.
We note that 
\begin{gather}
\label{ave}
\lim_{\ve\to 0} \tilde a_\ve(z)^i\ve^{-i(i+1)/2}=
\prod_{l=1}^i \tilde\Gamma_{p_l}(z), \\
\label{bve}
\lim_{\ve\to 0} \tilde b_\ve(z)^j\ve^{-j(j+1)/2}=
\prod_{r=1}^j \tilde\Gamma_{q_r}(z),
\end{gather} 
where $1\le i\le \la_0$ and $1\le j\le s$. 
Denote the right hand side of $(\ref{ave})$ by $\tilde a^{[i]}(z)$
and the right hand side of $(\ref{bve})$ by $\tilde b^{[j]}(z)$.
We have
$$\tilde a^{[i]}(z)=\tilde\Gamma_{p_1+\ldots +p_i}(z),\quad 
  \tilde b^{[j]}(z)=\tilde\Gamma_{q_1+\ldots +q_j}(z).$$ 
Because of
\begin{gather}
(p_1+\ldots + p_i, p_1+\ldots + p_j)=
(q_1+\ldots + q_i, q_1+\ldots + q_j)=2\min(i,j),\\
(p_1+\ldots + p_i, q_1+\ldots + q_j) =\max(0,i-\la_j)
\end{gather}
we obtain that the character of the algebra generated by 
$\tilde a^{[i]}(z)$, $\tilde b^{[j]}(z)$ is equal to the right hand side of
$(\ref{gmf})$ (see Proposition $\ref{charvo}$). 
This gives that $\ch_{z,u,q} A_{\bl;\bc,\bd}$ is greater than or
equal to the right hand side of $(\ref{gmf})$. 
To finish the proof we use Lemma $\ref{charlem}$ which gives the upper
bound for the character of $A_{\bl;\bc,\bd}$, which coincides with the right
hand side of $(\ref{gmf})$. Theorem is proved.
\end{proof}

\begin{cor}
The statement of Theorem $\ref{mth}$ is true. 
\end{cor}

We will need  the following  Corollary from the proof of the previous theorem. 
\begin{cor}
\label{vor}
Vertex operators $(\ref{point})$ provides a vertex operator realization of the 
algebra $A_{\la;\bc,\bd}$.
\end{cor}

\section{Fusion products}
\subsection{Principal subspaces}
In this section we apply results of the previous section to the study of 
fusion products of integrable irreducible representations of $\slth$.
We first recall main definitions.

Let $\g$ be some Lie algebra, $V_1,\ldots, V_n$ be its cyclic representations
with cyclic vectors $v_1,\ldots, v_n$, $Z=(z_1,\ldots, z_n)$ be the set
of pairwise distinct complex numbers. 
Denote by $V_i(z_i)$ the corresponding evaluation representations of 
$\g\otimes\C[u]$.
We consider a filtration $F_l$ on the tensor product 
$\bigotimes_{j=1}^n V_j(z_j)$ defined by
$$F_l=\mathrm{span}\{ x_1\otimes u^{i_1}\cdots x_p\otimes u^{i_p}
(v_1\T \ldots \T v_n):\
i_1+\ldots + i_p \le l, x_j\in\g\}.$$
The adjoint graded representation of $\g\otimes\C[u]$ is called the fusion 
product of $V_i$ and is denoted by
$V_1(z_1)*\ldots *V_n(z_n)$. Note that fusion products possess natural
$u$-grading defined by $\deg v=l$ for $v\in F_l/F_{l-1}$.

\begin{conj}
Let $\g$ be an affine Kac-Moody algebra, $V_i$ be its integrable 
representations. 
Then the corresponding fusion product does not depend on $Z$ as a 
representation  of $\g\otimes\C[u]$.
\end{conj}

This conjecture is obvious in the case $n=2$ (for an arbitrary Lie algebra 
$\g$). In what follows we study the case $\g=\slth$ and $n=2$.

We fix some notations first. Let 
$\slth=\slt\otimes\C[t,t^{-1}]\oplus \C K \oplus\C d$, where $K$ is the central 
element, and $d$ is the degree element. Let $e,h,f$ be the standard basis of 
$\slt$.
For $x\in\slt$ we set $x_i=x\T t^i\in\slth$.
Let $V_1=L_{i_1,k_1}$ and $V_2=L_{i_2,k_2}$ be two integrable irreducible 
$\slth$-modules with $0\le i_1\le k_1$, $0\le i_2\le k_2$. 
We fix a vector $v_{i,k}\in L_{i,k}$ which satisfy 
\begin{gather*}
f_0v_{i,k}=e_1v_{i,k}=0,\qquad U(\slth)v_{i,k}=L_{i,k},\\
h_0 v_{i,k}=-iv_{i,k},\ K v_{i,k}=kv_{i,k},\ dv_{i,k}=0,
\end{gather*}
($v_{i,k}$ is a highest weight vector with respect to the nilpotent algebra 
of annihilation operators generated by $f_0$ and $e_1$).
The principal subspace 
$W_{i,k}\hk L_{i,k}$ is defined by 
$$W_{i,k}=\C[e_0,e_{-1},\ldots]\cdot v_{i,k}.$$
This space is $z,q$ bi-graded by
$$\deg_z e_i=1,\ \deg_q e_i=-i,\ \deg_z v_{i,k}=\deg_q v_{i,k}=0.$$
$W_{i,k}$ is a representation of an abelian algebra $\A$ spanned by
$e_i$, $i\le 0$. So the fusion product $W_{i_1,k_1}* W_{i_2,k_2}$ is
a representation of $\A\oplus \A\T u$. Let 
$w_{i_1,k_1;i_2,k_2}\in W_{i_1,k_1}* W_{i_2,k_2}$ be the 
image of $v_{i_1,k_1}\T v_{i_2,k_2}$.
Then the action of the operators $e_i$ and $e_i\T u$ ($i\le 0$) on 
$w_{i_1,k_1;i_2,k_2}$ generates the fusion product.  Our goal is
to describe the ideal of relations in $\C[e_i,e_i\T u]_{i\le 0}$.
We will show that this ideal is isomorphic to $I_{\la;\bc,\bd}$ for 
special values of parameters. 

Namely let $I_{i_1,k_1;i_2,k_2}\hk \C[e_j, e_j\T u]_{j\le 0}$ be the ideal of 
relations, i.e.
$$W_{i_1,k_1}* W_{i_2,k_2}\simeq  \C[e_i, e_i\T u]_{i\le 0}/I_{i_1,k_1;i_2,k_2}.$$
\begin{prop}
\label{prisom}
An isomorphism  of algebras
\begin{equation}
\label{isom}
\C[e_j, e_j\T u]_{j\le 0}\to \C[a_j,b_j]_{i\le 0},\quad
e_i\mapsto a_i, \ e_i\T u\mapsto b_i
\end{equation}
induces the isomorphism of ideals
$$I_{i_1,k_1;i_2,k_2}\to I_{\la^{(k_1,k_2)};\delta^{(i_1+i_2+1)},
\delta^{(\min(i_1,i_2)+1)}},$$
where $\la^{(k_1,k_2)}=(k_1+k_2,k_1+k_2-2,\ldots, |k_1-k_2|)$ and 
$(\delta^{(i)})_j=\delta_{i,j}$.
\end{prop}

We start with the following lemma.

\begin{lem}
\label{surj}
We have an embedding of ideals (with respect to the isomorphism $(\ref{isom})$)
$$I_{i_1,k_1;i_2,k_2}\hk I_{\la^{(k_1,k_2)};\delta^{(i_1+i_2+1)},
\delta^{(\min(i_1,i_2)+1)}},$$
or the surjection
$$A_{\la^{(k_1,k_2)};\delta^{(i_1+i_2+1)},\delta^{(\min(i_1,i_2)+1)}}\to
W_{i_1,k_1}* W_{i_2,k_2}.$$
\end{lem}
\begin{proof}
Let $I_{i,k}\hk \C[e_j]_{j\le 0}$ be the ideal of relations in $W_{i,k}$,
i.e. 
$$I_{i,k}=\{p(e_0,e_{-1},\ldots): p v_{i,k}=0\}, \ 
W_{i,k}\simeq \C[e_j]_{j\le 0}/ I_{i,k}.$$ 
Denote $e(z)=\sum_{i\le 0} e_i z^{-i}$.
Then $I_{i,k}$ is generated by
coefficients of series $e(z)^{k+1}$ and first $l-i$ coefficients of the 
series $e(z)^l$, $l>i$ (the defining relations in $W_{i,k}$ are $e(z)^{k+1}=0$
and $e(z)^l\div z^{l-i}$, $l>i$).  
We recall that $W_{i_1,k_1}* W_{i_2,k_2}$ is
an adjoint graded space with respect to the filtration on the tensor product of 
evaluation representations $W_{i_1,k_2}(z_1)\T W_{i_2,k_2}(z_2)$, $z_1\ne z_2$.
We set
$$e(z)=\sum_{j\ge 0} e_{-j}z^j,\quad  
e^{(1)}(z)=e(z)\T \Id, \quad e^{(2)}(z)=\Id\T e(z).$$
Then obviously 
\begin{equation}
\label{i1+i2}
(e^{(1)}(z)+e^{(2)}(z))^l (v_{i_1,k_1}\T v_{i_2,k_2}) \div 
z^{l-i_1-i_2} \text{ for } l> i_1+i_2
\end{equation}
in the tensor product $W_{i_1,k_1}\T W_{i_2,k_2}$. 
This gives $(\ref{*})$ for $\bc=\delta^{(i_1+i_2+1)}$. Now let $z_1=1, z_2=0$.
Then $e_j\T u$ acts on the tensor product $W_{i_1,k_1}(1)\T W_{i_2,k_2}(0)$
as $e_j\T \Id$. We state that in  
$W_{i_1,k_2}* W_{i_2,k_2}$ the following equation is true:
\begin{equation}
\label{mini1i2}
e^{(1)}(z)^j\div z^{j-\min(i_1,i_2)} \text{ for } j> \min(i_1,i_2). 
\end{equation}
In fact, $e^{(1)}(z)^l\div z^{l-i_1}$. We also have
\begin{multline*}
e^{(2)}(z)^l=(e^{(1)}(z)+e^{(2)}(z)-e^{(1)}(z))^l=\\ (-e^{(1)}(z))^l+
\sum_{j=1}^l \genfrac{(}{)}{0pt}{0}{l}{j} 
(-e^{(1)}(z))^{l-j} (e^{(1)}(z)+e^{(2)}(z))^j.
\end{multline*}
Therefore  in $W_{i_1,k_2}* W_{i_2,k_2}$ we have
$$e^{(1)}(z)^l=(-e^{(2)}(z))^j\div z^{l-i_2}.$$
This gives $(\ref{**})$ for $\cd=\delta^{(\min(i_1,i_2)+1)}$.

We now check that
\begin{equation}
\label{lak1k2}
(e^{(1)}(z)+e^{(2)}(z))^{k_1+k_2-2i+1} e^{(1)}(z)^i=0,\ i=0,1,\ldots,
\min(k_1,k_2).
\end{equation}
Consider one-dimensional fusion product 
$(\C[e]/e^{k_1+1}) * (\C[e]/e^{k_2+1})$. As shown in \cite{FF1} the relations 
in this
fusion product are of the form
$$(e^{(1)}+e^{(2)})^{k_1+k_2-2i+1} (e^{(1)})^i=0,\ i=0,1,\ldots,
\min(k_1,k_2),$$
where $e^{(1)}=e\T\Id$, $e^{(2)}=\Id\T e$. This gives $(\ref{lak1k2})$.
Combining together  $(\ref{i1+i2})$, $(\ref{mini1i2})$ and $(\ref{lak1k2})$
we obtain our lemma. 
\end{proof}

Now our goal is to prove that 
$$
W_{i_1,k_1}* W_{i_2,k_2}\simeq
A_{\la^{(k_1,k_2)};\delta^{(i_1+i_2+1)},\delta^{(\min(i_1,i_2)+1)}}.
$$
We use the fermionic realization of $W_{i,k}$ (see \cite{FF2}). 
Let us briefly recall the construction.

Let $\phi_i, \psi_i$, $i\ge 0$ be anticommuting variables, 
$$\phi(z)=\sum_{i\ge 0} \phi_{-i} z^i,\quad 
\psi(z)=\sum_{i\ge 0} \psi_{-i} z^i.$$
We denote by $\Lambda=\Lambda(\phi_i,\psi_i)_{i\le 0}$ an exterior algebra
in variables $\phi_i, \psi_i$. 
Then there exists an
embedding $\imath_{i,k} :W_{i,k}\hk \Lambda^{\T k}$ such that
\begin{equation}
\label{vik}
v_{i,k}\mapsto \wt v_{i,k}=
\underbrace{1 \T\ldots\T 1}_i \T \phi_0\T\ldots\T \phi_0.
\end{equation}
This embedding is defined via the identification 
$e(z)=\sum_{n=1}^k \stackrel{n}{\phi}(z)\stackrel{n}{\psi}(z)$,  where
$$\stackrel{n}{\phi}(z)=\underbrace{\Id\T\ldots\T\Id}_{n-1}\T 
\phi(z)\T\Id\T\ldots\T\Id$$
and similarly for $\stackrel{n}{\psi}(z)$. 
Note that 
$$\imath_{i,k} (W_{i,k})\hk \Lambda_{ev}\cdot \wt v_{i,k},$$
where $\Lambda_{ev}$ is an even part of $\Lambda$ generated by the products
$\psi_i\psi_j$, $\psi_i\phi_j$, $\phi_i\phi_j$.
The algebra $\Lambda_{ev}$ is naturally $(z,q)$ bi-graded with
$$\deg_z \psi_i\psi_j=1,\    \deg_q \psi_i\psi_j=i+j$$
and similarly for other generators. Fixing 
$\deg_z \wt v_{i,k}= \deg_q \wt v_{i,k}=0$ we obtain a bi-grading on 
$\Lambda_{ev}\cdot \wt v_{i,k}$ such that $\imath_{i,k}$ is an
embedding of $(z,q)$ bi-graded spaces.

Note that our construction gives an embedding
$$\imath_{i_1,k_1}\T \imath_{i_2,k_2}: 
W_{i_1,k_1}\T W_{i_2,k_2}\hk \Lambda^{\T (k_1+k_2)}.$$
Our goal is to construct a continuous family of algebras 
$B(\ve)\hk \Lambda^{\T (k_1+k_2)}$, 
$0\le \ve \le 1$ which "connects" $W_{i_1,k_1}\T W_{i_2,k_2}$ and 
$A_{\la^{(k_1,k_2)};\delta^{(i_1+i_2+1)},\delta^{(\min(i_1,i_2)+1)}}$ inside
$\Lambda^{\T (k_1+k_2)}$. So we need a fermionic realization of an algebra
$$A_{\la^{(k_1,k_2)};\delta^{(i_1+i_2+1)},\delta^{(\min(i_1,i_2)+1)}}.$$

\begin{lem}
The map $\jmath$ defined by 
\begin{gather} 
\label{1to}
1\mapsto \wt v_{i_1,k_2}\T \wt v_{i_2,k_2},\\ \label{a(z)}
a(z)\mapsto \wt a(z)=
\sum_{n=1}^{k_1+k_2} \stackrel{n}{\phi}(z)\stackrel{n}{\psi}(z),\\ \label{b(z)}
b(z)\mapsto \wt b(z)=
\begin{cases}
\sum_{n=1}^{k_1} \stackrel{n}{\phi}(z)\stackrel{n+k_1}{\psi}(z), \ i_1\le i_2,\\
\sum_{n=1}^{k_1} \stackrel{n+k_1}{\phi}(z)\stackrel{n}{\psi}(z), \ i_1>i_2.
\end{cases}
\end{gather}
provides an embedding
$$A_{\la^{(k_1,k_2)};\delta^{(i_1+i_2+1)},\delta^{(\min(i_1,i_2)+1)}}\hk
\Lambda^{\T (k_1+k_2)}.$$
\end{lem}
\begin{proof}
Our Lemma is an immediate consequence of vertex operator realization of 
$A_{\la^{(k_1,k_2)};\delta^{(i_1+i_2+1)},\delta^{(\min(i_1,i_2)+1)}}$ (see
Corollary $\ref{vor}$) and a version of the boson-fermion correspondence 
(see for example \cite{BF}). We give a sketch of the proof here.

Recall that the boson-fermion correspondence allows to realize lattice
VOAs inside the space build up from the fermionic particles of the type
$\phi_i$ and $\psi_i$. In particular the space of states of corresponding
VOA (the direct sum of Fock modules) is replaced by the space of semi-infinite
forms or by the tensor products of such spaces.
In this realization fields $\phi(z)$ and $\phi(z)\psi(z)$ 
correspond to the one-dimensional vertex operators $\Gamma_\al(z)$ with
the length of $\al$  equals $1$ or $2$ respectively. To proceed to 
multi-dimensional even case one must consider independent fileds
of the type $\stackrel{n}{\phi}(z)\stackrel{m}{\psi}(z)$. 
The independence means that these field represent
vertex operators $\Gamma_{\be_{n,m}}(z)$ and the scalar products are given by
$(\be_{n,m},\be_{k,l})=\delta_{n,k}+\delta_{m,l}.$

We now apply this boson-fermion correspondence to the proof of our lemma. 
Let $p_1,\ldots, p_{k_1+k_2}$, $q_1,\ldots,q_{k_1}$ be vectors such that
\begin{gather} 
\stackrel{n}{\phi}(z)\stackrel{n}{\psi}(z)=\Gamma_{p_{k_1+k_2+1-n}} (z),
1\le n\le k_1+k_2\\
\stackrel{m}{\phi}(z)\stackrel{m+k_1}{\psi}(z)=\Gamma_{q_m}(z), \ i_1\le i_2,
\quad \stackrel{m+k_1}{\phi}(z)\stackrel{m}{\psi}(z)=\Gamma_{q_m}(z), \ 
i_1>i_2,
\end{gather}
where $1\le m\le k_1$.
Then $(p_n,q_m)=0$  unless $k_1+k_2-2(m-1)\ge n > k_1+k_2-2m$. In the latter
case the corresponding scalar product equals $1$. Therefore according to 
Corollary $\ref{vor}$ to finish the proof of our Lemma we only need to check 
the initial conditions
$$\tilde a(z)^{i_1+i_2+l} \tilde v_{i_1,k_1}\T  \tilde v_{i_2,k_2} \div z^l,\quad
\tilde b(z)^{\min(i_1,i_2)+l} \tilde v_{i_1,k_1}\T  \tilde v_{i_2,k_2} 
\div z^l.$$
But this just follows from $\stackrel{n}{\phi}(z) \tilde v_{i,k}\div z^l$
for $n>i$ (see $(\ref{vik})$).
\end{proof}

{\bf Proof of Proposition \ref{prisom}: }
\begin{equation}
\label{isomfus}
W_{i_1,k_1}* W_{i_2,k_2}\simeq
A_{\la^{(k_1,k_2)};\delta^{(i_1+i_2+1)},\delta^{(\min(i_1,i_2)+1)}}.
\end{equation}
\begin{proof}
Because of Lemma $\ref{surj}$ it suffices to check that the character of the
left hand side of $(\ref{isomfus})$ coincides with the 
character of the right hand side.  

We construct a continuous family of algebras 
$B(\ve)\hk \Lambda^{\T (k_1+k_2)}$, 
$0\le \ve <1$ which "connects" $W_{i_1,k_1}\T W_{i_2,k_2}$ and 
$A_{\la^{(k_1,k_2)};\delta^{(i_1+i_2+1)},\delta^{(\min(i_1,i_2)+1)}}$.
We want $B(\ve)$ to satisfy
\begin{description} 
\item[A] $B(0)=(\imath_{i_1,k_1}\T \imath_{i_2,k_2}) (W_{i_1,k_1}\T W_{i_2,k_2})$
\item[B] $B(\ve)\simeq B(0)$ as bi-graded vector spaces
\item[C] $\lim_{\ve\to 1} B(\ve)=\jmath  
A_{\la^{(k_1,k_2)};\delta^{(i_1+i_2+1)},\delta^{(\min(i_1,i_2)+1)}}$
\end{description}
Note that the existence of such deformation proves our proposition.

Set
\begin{equation}
\label{b}
\wt b^\ve(z)=
\begin{cases}
\sum_{n=1}^{k_1} \stackrel{n}{\phi}(z)
(\ve\stackrel{n}{\psi}(z)+(1-\ve)\stackrel{n+k_1}{\psi}(z)), \ i_1\le i_2,\\
\sum_{n=1}^{k_1} 
(\ve\stackrel{n}{\phi}(z)+(1-\ve)\stackrel{n+k_1}{\phi}(z)) 
\stackrel{n}{\psi}(z), \ i_1>i_2.
\end{cases}
\end{equation}
We denote by $B(\ve)\hk \Lambda^{\T (k_1+k_2)}$ the subspace generated by
the coefficients of $\wt a(z)$ and $b^\ve(z)$ from the vector 
$\wt v_{i_1,k_1}\T \wt v_{i_2,k_2}$. We now check ${\bf A}$, ${\bf B}$ and 
${\bf C}$.

${\bf A}$ is obvious because $\wt b^0(z)$ reduces to the current $\wt b(z)$. 

We now check ${\bf B}$. 
First note that for any $0< \ve_1\le \ve_2<1$ we have $B(\ve_1)\simeq 
B(\ve_2)$. In fact, fix some $0<\ve <1$ and redefine
\begin{gather*}
1/2\stackrel{n}{\psi'}(z)=(1-\ve)\stackrel{n}{\psi}(z),\ 1\le n\le k_1,\\
1/2\stackrel{m}{\psi'}(z)=\ve\stackrel{m}{\psi}(z),\ k_1+1\le m\le k_1+k_2,\\
2\stackrel{n}{\phi'}(z)=(1-\ve)^{-1}\stackrel{n}{\phi}(z),\ 1\le n\le k_1,\\
2\stackrel{m}{\phi'}(z)=\ve^{-1}\stackrel{m}{\phi}(z), \ k_1+1\le m\le k_1+k_2,
\end{gather*}
for $i_1\le i_2$ and 
\begin{gather*}
1/2\stackrel{n}{\phi'}(z)=(1-\ve)\stackrel{n}{\phi}(z),\ 1\le n\le k_1,\\
1/2\stackrel{m}{\phi'}(z)=\ve\stackrel{m}{\phi}(z),\ k_1+1\le m\le k_1+k_2, \\
2\stackrel{n}{\psi'}(z)=(1-\ve)^{-1}\stackrel{n}{\psi}(z),\ 1\le n\le k_1,\\
2\stackrel{m}{\psi'}(z)=\ve^{-1}\stackrel{m}{\psi}(z), \ k_1+1\le m\le k_1+k_2,
\end{gather*}
for $i_1>i_2$. Then one can easily check that $\wt a(z)$ doesn't change and 
$\wt b^\ve(z)$ becomes
$\wt b^{1/2}(z)$ (up to the nonzero constant). Therefore, $B(\ve)\simeq B(1/2)$
for any $0<\ve <1$.  

Now note that 
\begin{equation}
\label{sep}
(2\wt b^{1/2} (z))^{k_1+1} \wt v_{i_1,k_1}\T \wt v_{i_2,k_2}=
(a(z)-2\wt b^{1/2} (z))^{k_2+1} \wt v_{i_1,k_1}\T \wt v_{i_2,k_2}=0
\end{equation} 
and 
\begin{gather}
\label{i1}
(2\wt b^{1/2} (z))^l \wt v_{i_1,k_1}\T \wt v_{i_2,k_2}\div z^{l-i_1+1}, \ l>i_1,\\
\label{i2}
(\wt a(z) - 2\wt b^{1/2} (z))^l \wt v_{i_1,k_1}\T \wt v_{i_2,k_2}\div 
z^{l-i_2+1}, \ l>i_2.
\end{gather}
In fact, one has 
$$2\wt b^{1/2}(z)=
\begin{cases}
\sum_{n=1}^{k_1} \stackrel{n}{\phi}(z)
(\stackrel{n}{\psi}(z)+\stackrel{n+k_1}{\psi}(z)), \ i_1\le i_2,\\
\sum_{n=1}^{k_1} 
(\stackrel{n}{\phi}(z)+\stackrel{n+k_1}{\phi}(z)) 
\stackrel{n}{\psi}(z), \ i_1>i_2;
\end{cases}, $$
and
$$\wt a(z)- 2\wt b^{1/2}(z)=
\begin{cases}
\sum_{n=k_1+1}^{k_1+k_2} (\stackrel{n}{\phi}(z)- \stackrel{n-k_1}{\phi}(z))
\stackrel{n}{\psi}(z), \ i_1\le i_2,\\
\sum_{n=k_1+1}^{k_1+k_2} 
\stackrel{n}{\phi}(z) 
(\stackrel{n}{\psi}(z)- \stackrel{n-k_1}{\psi}(z)), \ i_1>i_2.
\end{cases} 
$$
Now $(\ref{sep}), (\ref{i1}), (\ref{i2})$ follows from the above exact 
expressions of the currents. Relations 
$(\ref{sep}), (\ref{i1}), (\ref{i2})$  provide a surjection $B(0)\to B(1/2)$
($2\wt b^{1/2}(z)$ and $\wt a(z)$ correspond to the currents
$e^{(1)}(z)$ and $e^{(2)}(z)$ in $B(0)\simeq W_{i_1,k_1}\T W_{i_2,k_2}$). 
This gives ${\bf B}$ (because
all $B(\ve)$ with $0<\ve<1$ are isomorphic and our deformation is
continuous).

To show ${\bf C}$ we note that in view of formulas $(\ref{1to}), (\ref{a(z)}), 
(\ref{b(z)})$ and $(\ref{b})$ we have an embedding
$$\jmath A_{\la^{(k_1,k_2)};\delta^{(i_1+i_2+1)},\delta^{(\min(i_1,i_2)+1)}}\hk
   \lim_{\ve\to 1} B(\ve).$$
But all $B(\ve)$ (including $B(0)$) are isomorphic and 
$$\ch_{z,q} B(0)\le 
\ch_{z,q} A_{\la^{(k_1,k_2)};\delta^{(i_1+i_2+1)},\delta^{(\min(i_1,i_2)+1)}}.$$
Therefore, using 
Lemma $\ref{surj}$ we obtain ${\bf C}$.
\end{proof}

\begin{cor}
\label{w1w2}
Let $k_1\le k_2$. Then
\begin{multline}
\ch_{z,u,q} (W_{i_1,k_1}* W_{i_2,k_2})=
\sum_{\bn\in \Z^{k_1+k_2}_{\ge 0}, \bm\in\Z_{\ge 0}^{k_1}}
u^{|\bm|} z^{|\bn|+|\bm|}\times\\ 
q^{-|\bn|-|\bm|+\sum\limits_{j=i_1+i_2+1}^{k_1+k_2} (j-i_1-i_2)n_j +
\sum\limits_{j=\min(i_1,i_2)+1}^{k_1} (j-\min(i_1,i_2))m_j}\times\\
\frac{q^{\bn A\bn/2 + \bn B\bm + \bm A\bm/2}}
{(q)_\bn (q)_\bm},
\end{multline}
where $A_{i,j}=2\min(i,j)$ and $B_{i,j}=\max(0,i-k_1-k_2+2j)$.
\end{cor}

\subsection{The limit construction.} 
In this subsection we derive a fermionic formula 
for the character of the fusion product $L_{i_1,k_1}* L_{i_2,k_2}$. 
Note that $L_{i,k}$ is bi-graded by operators $d$ and $h_0$. Therefore
for the fusion product the $z,u,q$ character is naturally defined.

Let  $v_{i,k}^N\in L_{i,k}$, $N\in\Z$ be the set of extremal vectors
(the weight of $v_{i,k}^N$ is a weight of $v_{i,k}$ shifted by the $N$-th
power of the translation element from the Weyl group of $\slth$). We fix
$h_0v^N_{i,k}=(-i-2Nk)v_{i,k}^N$. 
Introduce the $N$-th principal subspace
$W_{i,k}^N\hk L_{i,k}$ by
$$W_{i,k}^N=\C[e_{2N}, e_{2N-1},\ldots]\cdot v_{i,k}^N$$
(note that $e_{2N+1}v_{i,k}^N=0$).
We recall that there exists an isomorphism
\begin{equation}
\label{shift}
W_{i,k}\simeq W_{i,k}^N,\quad v_{i,k}\mapsto v_{i,k}^N,
e_j\to e_{j+2N}
\end{equation}
and
$$L_{i,k}=W_{i,k}^0\hk W_{i,k}^1 \hk W_{i,k}^2 \hk\ldots.$$
Using this limit construction we can write a fermionic formula for the 
character of the fusion product of integrable modules.

\begin{cor}
\begin{multline}
\label{limform}
\ch_{z,u,q} (L_{i_1,k_1}* L_{i_2,k_2})= \lim_{N\to\infty}
z^{-i_1-i_2-2N(k_1+k_2)} q^{N^2(k_1+k_2)+N(i_1+i_2)}\times\\
\sum_{\bn\in \Z^{k_1+k_2}_{\ge 0}, \bm\in\Z_{\ge 0}^{k_1}}
u^{|\bm|} z^{2(|\bn|+|\bm|)} q^{(-2N-1)(|\bn|+|\bm|)}\times\\
q^{\sum\limits_{j=i_1+i_2+1}^{k_1+k_2} (j-i_1-i_2)n_j +
\sum\limits_{j=\min(i_1,i_2)+1}^{k_1} (j-\min(i_1,i_2))m_j}
\frac{q^{\bn A\bn/2 + \bn B\bm + \bm A\bm/2}}
{(q)_\bn (q)_\bm}.
\end{multline}
\end{cor}
\begin{proof}
Follows from  Corollary \ref{w1w2}, isomorphism $(\ref{shift})$  and equalities
$$h_0 v_{i,k}^N=(-i-2Nk)v_{i,k}^N,\quad dv_{i,k}^N=(kN^2+Ni)v_{i,k}^N.$$
We only note that a 
factor $2$ in $z^{2(|\bn|+|\bm|)}$ comes from the relation $[h_0,e_i]=2e_i$
(in $W_{i,k}$ we had $\deg_z e_i=1$) and the power
$q^{(-2N-1)(|\bn|+|\bm|)}$ comes from the shift $e_j\to e_{j+2N}$ (see 
$(\ref{shift})$).
\end{proof}

Introduce a new set of variables $s_i$ instead of $n_i$ in 
$(\ref{limform})$ by
$$n_i=s_i-s_{i+1}, 1\le i <k_1+k_2,\ n_{k_1+k_2}=s_{k_1+k_2}+N-
\frac{|\bm|}{k_1+k_2}.$$
Then for the power of $z$ in $(\ref{limform})$ we have
$$-i_1-i_2-2N(k_1+k_2)+2(|\bn|+|\bm|)=-i_1-i_2+2\sum_{i=1}^{k_1+k_2} s_i.$$
We now rewrite the power of $q$ in $(\ref{limform})$ in terms of new variables.
Note that 
\begin{gather*}
\bn A\bn/2=\sum_{i=1}^{k_1+k_2} (n_i+\ldots +n_{k_1+k_2})^2=
\sum_{i=1}^{k_1+k_2} (s_i+N-\frac{|\bm|}{k_1+k_2})^2,\\
\bn B\bm=|\bm|(2N-2\frac{|\bm|}{k_1+k_2})+\sum_{i+2j\ge k_1+k_2+1} m_j s_i,\\
|\bn|+|\bm|=N(k_1+k_2)+\sum_{i=1}^{k_1+k_2} s_i.
\end{gather*}
Therefore the power of $q$ in $(\ref{limform})$ is equal to
\begin{multline}
\label{P}
\sum_{i=1}^{k_1+k_2} s_i^2-
\frac{|\bm|}{k_1+k_2}(|\bm|+k_1+k_2-i_1-i_2+ 2\sum_{i=1}^{k_1+k_2} s_i)+\\
\sum_{i+2j\ge k_1+k_2+1} m_j s_i+  \bm A\bm/2 - \sum_{i=1}^{i_1+i_2} s_i +
\sum_{j=\min(i_1,i_2)+1}^{k_1} (j-\min(i_1,i_2))m_j.
\end{multline}
(Note that the powers which contain $N$ cancel.) 
We thus obtain the following theorem
\begin{theorem}
\begin{multline*}
\ch_{z,u,q} L_{i_1,k_1}* L_{i_2,k_2}=\\  \frac{1}{(q)_\infty}
\sum_{\substack{\bs\in\Z^{k_1+k_2}, \bm\in\Z_{\ge 0}^{k_1}\\ 
s_1\ge \ldots \ge s_{k_1+k_2}}} z^{-i_1-i_2+2\sum_{i=1}^{k_1+k_2} s_i}
\frac{q^{P(\bs,\bm)}}{(q)_\bm \prod_{i=1}^{k_1+k_2-1} (q)_{s_i}} ,
\end{multline*}
where $P(\bs,\bm)$ is given by $(\ref{P})$ and 
$(q)_\infty=\prod_{i=1}^\infty (1-q^i)$.
\end{theorem}
\begin{proof}
We only note that the factor $(q)_\infty$ comes from the limit
$$\lim_{N\to\infty} (q)_{n_{k_1+k_2}}=
\lim_{N\to\infty}(q)_{s_{k_1+k_2}+N-\frac{|\bm|}{{k_1+k_2}}}.$$ 
\end{proof}

\newcounter{a}
\setcounter{a}{1}

\end{document}